\documentstyle[11pt,dgspp,epsf,psrotate]{article}

\newtheorem{theorem}{Theorem}[section]
\newtheorem{proposition}[theorem]{Proposition}
\newtheorem{lemma}[theorem]{Lemma}
\newtheorem{corollary}[theorem]{Corollary}
\newtheorem{definition}[theorem]{Definition}
\newremark{example}[theorem]{Example}
\newremark{remark}[theorem]{Remark}

\newcommand{\av}{a_{\scriptscriptstyle V}}
\newcommand{\bdel}[1]{{\bar\nabla}_{\!#1}}
\newcommand{\bU}{\bar{U}}
\newcommand{\bV}{\bar{V}}
\newcommand{\cD}{\cal D}
\newcommand{\cJ}{\escr J}
\newcommand{\cT}{\cal T}
\newcommand{\dc}{\dot{c}}
\newcommand{\ddc}{\ddot{c}}
\newcommand{\dx}{\dot{x}}

\newcommand{\ddx}{\ddot{x}}

\newcommand{\del}[1]{\nabla_{\!#1}}
\newcommand{\ds}{\oplus} 
\newcommand{\R}{{\Bbb R}}

\newcommand{\D}{\nabla}

\newcommand{\g}{\gamma}
\newcommand{\gl}{\frak{gl}}
\newcommand{\half}{\mbox{$\textstyle\frac{1}{2}$}}

\newcommand{\lsp}{[\kern-0.15em[} 
\newcommand{\nd}{\vdash}
\newcommand{\p}{\mbox{$\pi$}}
\newcommand{\rsp}{]\kern-0.15em]} 
\newcommand{\surj}{\rightarrow\kern-.82em\rightarrow}

\newcommand{\V}{{\rscr V}}
\newcommand{\ve}{\varepsilon}

\renewcommand{\H}{{\rscr H}}
\renewcommand{\P}{\Phi}
\renewcommand{\S}{\lower .3ex\hbox{$\escr S$}}

\makeatletter
\newcommand{\Ad}[1]{\mathop{\operator@font Ad}\nolimits_{#1}}
\newcommand{\Aut}{\mathop{\operator@font Aut}\nolimits}
\newcommand{\con}{\mathop{\operator@font con}\nolimits}
\newcommand{\diag}{\mathop{\operator@font diag}\nolimits}
\newcommand{\dom}{\mathop{\operator@font dom}\nolimits}
\newcommand{\End}{\mathop{\operator@font End}\nolimits}
\newcommand{\Hom}{\mathop{\operator@font Hom}\nolimits}
\newcommand{\im}{\mathop{\operator@font im}\nolimits}
\newcommand{\Op}{\mathop{\operator@font Op}\nolimits}
\newcommand{\Ric}{\mathop{\operator@font Ric}\nolimits}
\newcommand{\spr}{\mathop{\operator@font Spray}\nolimits}
\newcommand{\tr}{\mathop{\operator@font tr}\nolimits}
\makeatother

\preprint{DRP3}
\title{Generalized Sprays and Nonlinear Connections}
\author{L. Del Riego\thanks{Partially supported by CONACYT grant 
        26594-E.}}
\address{Facultad de Ciencias\\
   Zona Universitaria\\
   Universidad Aut\'onoma de San Luis Potos\'\i\\
   San Luis Potos\'\i, SLP\\
   78290 MEXICO\\
   lilia@galia.fc.uaslp.mx}
\author{Phillip. E. Parker}
\address{Mathematics Department\\
   Wichita State University\\
   Wichita KS 67260-0033\\
   USA\\
   phil@math.twsu.edu}
\date{4 April 2003}
\abstract{
The main purposes of this article are to extend our previous results on
homogeneous sprays \cite{DRP1} to arbitrary (generalized) sprays, to show
that locally diffeomorphic exponential maps can be defined for any
(generalized) spray, and to give a (possibly nonlinear) covariant
derivative for any (possibly nonlinear) connection.  In the process, we
introduce {\em vertically homogeneous\/} connections.  Unlike homogeneous
connections, these allow us to include Finsler spaces among the
applications.

We provide significant support for the prospect of studying nonlinear
connections {\em via\/} (generalized) sprays. One of the most important 
is our generalized APS correspondence.}
\msc{53C15}{53C22, 58E10}

\begin{document}

\maketitle


\section{Introduction}
An important class of systems of second order differential equations can be
represented as (generalized) sprays on a manifold $M$ with tangent bundle 
$TM \surj M$.  So far only quadratic sprays are well understood, and they
correspond to linear connections.  But nonlinear connections are of real
interest, especially in some newer applications \cite{AA,AIM,A,R,V}.

In Riemannian geometry, the (usual) geodesic spray, whose integral curves 
are the geodesics of the Levi-Civita connection, has played an important 
r\^ole; see, for example, \cite{BC,B}. In Finsler geometry, four main 
connections have been used, none of them linear: those of Cartan, Berwald,
Hasiguchi, and Chern \cite{BCS}.  Riemannian geometry has been a main
thread of mathematics over the last century, and Finsler geometry has
recently undergone a great revival. Applications of it now include 
modeling the singular sets of Monge-Amp\'ere equations \cite{R}, studying 
the manifold of Hamiltonians \cite{BPo,Po}, and modeling river flows and 
mountain slopes \cite{AIM}.

One of our motivations for this work was the desire to make a 
comprehensive theory of sprays and nonlinear connections which would 
include Riemannian and Finsler spaces as examples. We have recovered 
enough of the Riemannian results to be assured of the correctness of our 
approach; comparing the results for Finsler spaces of our methods with 
those of other methods will be the subject of future study.

Section \ref{rev} contains our notation, conventions, and a summary of our
earlier article \cite{DRP1}.  In Section \ref{exp} we present the new
exponential maps defined by (generalized) sprays.  Section \ref{cs}
describes some relations between (possibly nonlinear) connections and
(generalized) sprays and the associated (possibly nonlinear) covariant
derivatives and geodesics.  Section \ref{gcs} begins with the extension of
the main results of \cite{BP6} to (generalized) sprays, using our new
construction of (generalized) exponential maps.  It also includes the
extension of the main stability result of \cite{BP4,DRP1} to all
(generalized) sprays.

Throughout, all manifolds are smooth (meaning $C^{\infty}$), connected,
paracompact, and Hausdorff.

The authors thank CONACYT and FAI for travel and support grants, Wichita
State University and Universidad Aut\'onoma de San Luis Potos\'{\i} for
hospitality during the progress of this work, and J. Hebda and A. Helfer
for helpful conversations.  Del Riego also thanks M. Mezzino for writing a
Mathematica package for her use.

\section{Review and definitions}\label{rev}
A (general) {\em spray\/} on a manifold $M$ is defined as a projectable
section of the second-order tangent bundle $TTM \surj TM$.  This is
precisely the condition needed to define a second-order differential
equation \cite{BC,B}.  Recall that an integral curve of a vector field on
$TM$ is the canonical lift of its projection if and only if the vector field
is projectable.  For any curve $c$ in $M$ with tangent vector field $\dot
c$, this $\dot c$ is the canonical lift of $c$ to $TM$ and $\ddot c$ is the
canonical lift of $\dot c$ to $TTM$.  Then each projectable vector field $S$
on $TM$ determines a second-order differential equation on $M$ by $\ddot c =
S\circ \dot c$, and any such curve with $\dot c(s_{0}) = v_{0}\in
T_{c(s_{0})}M$ is a solution with initial condition $v_0$.  Solutions are
preserved under translations of parameter, they exist for all initial
conditions by the Cauchy theorem, and, as our manifolds are assumed to be
Hausdorff, each solution will be unique provided we take it to have maximal
domain; {\em i.e.,} to be inextendible \cite{BC,DRD,HS}.

Let $J$ be the canonical involution on $TTM$ and $C$ the Euler (or
Liouville) vector field.  We recall that in local coordinates, $J(x,y,X,Y) =
(x,X,y,Y)$ and $C:(x,y)\mapsto (x,y,0,y)$.  
\begin{definition}\label{dfs}
A section $S$ of $TTM$ over $TM$ is a {\em spray\/} when $JS = S$; that
is, when it can be expressed locally as $S:(x,y)\mapsto (x,y,y,\S(x,y))$.
\end{definition}
Before commenting on this definition, we must briefly digress to consider
the notion of homogeneity for functions.

Consider the equation $f(ax) = a^m f(x)$.  In projective geometry, for 
example, one usually requires this to hold only for $a\ne 0$.  We shall 
call this {\em homogeneous\/} of degree $m$.  In other areas, such as 
Euler's Theorem in analysis, one further restricts to $a>0$.  We shall 
call this {\em positively homogeneous\/} of degree $m$.  Finally, in order
that homogeneity of degree 1 coincide with linearity, one must allow any 
scalar $a\in\R$ (including zero).  We shall call this {\em completely 
homogeneous\/} of degree $m$ and denote it by $h(m)$.

The difference between homogeneity and complete homogeneity is minor; 
essentially, it is just the difference between working on $TM-0$ and on 
$TM$.  The difference between positive homogeneity and the other two is 
more significant.  For example, the inward-going and outward-going radial 
geodesics of the Finsler-Poincar\'e plane in \cite{BCS} have different 
arclengths.

Now we are ready to consider homogeneity for sprays.
\begin{definition}\label{hms}
We say that a spray $S$ is\/ {\em homogeneous} of degree $m$ when the
functions $\S(x,y)$ are completely homogeneous (respectively, homogeneous)
of degree $m$ in the vertical component in some induced local 
coordinates:
$\S(x,ay) = a^{m}\S(x,y)$ for some $m\ge 2$ (respectively, $m<2$) and all
scalars $a\in\R$ (respectively, $a\ne 0$).
\end{definition}
The break comes at $m=2$ because an $h(m)$ spray is to be associated with 
a connection whose homogeneity formula will contain $a^{m-2}$; {\em cf.} 
(\ref{vhceq}).  In the distinguished induced local coordinates,
$S:(x,ay)\mapsto (x,ay,ay,a^m \S(x,y))$.  Only induced local
coordinates $(x',y',X',Y')$ related to this $(x,y,X,Y)$ by a
block-diagonal transition matrix
$$ \left[ \begin{array}{cc} 
\frac{\partial x'}{\partial x} & \frac{\partial x'}{\partial y}\\[1ex]
\frac{\partial y'}{\partial x} & \frac{\partial y'}{\partial y} 
\end{array} \right] = 
\left[ \begin{array}{cc} 
\frac{\partial x'}{\partial x} & 0\\[1ex]
0 & \frac{\partial y'}{\partial y} 
\end{array} \right] $$
will preserve the form of such an $S$. Other induced local coordinates
preserve the correct degree of homogeneity in the vertical component $Y$,
but may change the degree of homogeneity in the ``horizontal" component
$X$.  Thus from now on, we shall use only these {\em admissible\/} atlases
on $TM$ when studying homogeneous sprays {\em et relata;} {\em cf.}
after Theorem \ref{cce}.
\begin{remark}
In the extant literature \cite{DRP1,G,G1,KV,LR}, one finds homogeneous
{\em vector fields\/} of degree $m$ defined by $[C,S] = (m-1)S$.  In {\em
any\/} (not just admissible) local coordinates, $S:(x,ay)\mapsto
(x,ay,a^{m-1}y,a^m \S(x,y))$.  It follows that a spray in our theory can
be a homogeneous vector field only for $m=2$.
\end{remark}
Hereinafter we shall call $h(2)$ sprays {\em quadratic\/} sprays, in
agreement with \cite{G1,KV,LR}.  (Note that {\em complete\/} homogeneity
is required for our quadratic sprays to coincide with the usual spray of
\cite{APS}.)  We denote the set of our sprays on $M$ by ${\spr}(M)$ and
those which are $h(m)$ by ${\spr}_m(M)$.  It has been usual to consider
only (positive) integral degrees of homogeneity, but we make no such
restriction.

Previously \cite{KV}, projectable vector fields on $TM$ were called
{\em semi\-sprays\/} and the name {\em sprays\/} used for those that were
homogeneous.  We will associate one of our (general) sprays to each
(possibly nonlinear) connection as its {\em geodesic spray\/} (see
Theorems \ref{cis} and \ref{cg=sg}), so we are using the name
``sprays" to reflect this new, extended r\^ole.  We do, however,
explicitly consider only sprays defined on the entire tangent bundle $TM$;
others \cite{A,BCS,KV} have used the reduced tangent bundle with the
0-section removed, which is appropriate when considering $h(m)$ sprays
when $m<2$ (including $m<0$).  For $0 \le m < 2$, one usually requires
sprays to be $C^0$ across the zero-section; {\em e.\,g.,} for Finsler
spaces.  Most of our results are easily seen to hold {\em mutatis
mutandis\/} in these cases as well; any unobvious exceptions will be noted
specifically.

In fact, the desire to make a comprehensive theory including Finsler
spaces was one of our motivations.  What {\em should\/} be the
Finsler-geodesic spray associated with a Finsler metric tensor is {\em
not\/} a homogeneous vector field, but an $h(1)$ spray in our theory; {\em
cf.} \cite{DR1} for related results.  However, what is frequently used as
the Finsler-geodesic spray has both quadratic and $h(1)$ parts; {\em cf.}
\cite{BCS}, for example.  We plan to address these peculiarities of the
existing theory in subsequent work.

Several important results concerning quadratic sprays \cite{APS,BC,D,KV}
rely on the facts that each such spray $S$ determines a unique torsion-free
linear connection $\Gamma$, and conversely, every quadratic spray $S$ arises
from a linear connection $\Gamma$ the torsion of which can be assigned
arbitrarily.  The solution curves of the differential equation $\ddot c =
S_{\Gamma}\circ \dot c$ for a connection-induced spray are precisely the
geodesics of that (linear) connection.  These solution curves are not only
preserved under translations, as is true in general, but also under affine
transformations of the parameter $s \mapsto as+b $ for constants $a,b$ with
$a\neq 0$.  Note that, with our definition, the latter also holds
for homogeneous sprays.  

In the general case, a (possibly nonlinear) connection $\Gamma$ gives rise
to a spray $S$ (see Proposition \ref{cis}), but the correspondence has not
been well studied before.  We shall extend most of the preceding features
of the quadratic spray---linear connection correspondence to the general
setting.  One of our ultimate goals is to determine just how well
nonlinear connections can be studied {\em via\/} sprays.

We continue with the principal definitions. Let $S$ be a (generalized)
spray on $M$.
\begin{definition}
We say that a curve $c:(a,b)\rightarrow M$ is a\/ {\em geodesic} of $S$ or 
an {\em $S$-geodesic} if and only if the natural lifting $\dot c$ of $c$ to
$TM$ is an integral curve of $S$.
\end{definition}
This means that if $\ddot c$ is the natural lifting of $\dot c$ to $TTM$, 
then $\ddot c = S(\dot c)$. 
\begin{definition}
We say that $S$ is\/ {\em pseudoconvex} if and only if for each compact $K
\subseteq M$ there exists a compact $K^\prime \subseteq M$ such that each
$S$-geodesic segment with both endpoints in $K$ lies entirely within
$K^{\prime}$.
\end{definition}
If we wish to work directly with the integral curves of $S$, we merely
replace ``in''  and ``within'' by ``over''.
\begin{definition} 
We say that $S$ is\/ {\em disprisoning} if and only if no inextendible
$S$-geodesic is contained in (or lies over) a compact set of $M$.
\end{definition}
In relativity theory, such inextendible geodesics are said to be imprisoned
in compact sets; hence our name for the negation of this property.

Following this definition, we make a convention:  all $S$-geodesics are 
always to be regarded as extended to the maximal parameter intervals ({\em
i.e.,} to be inextendible) unless specifically noted otherwise.  When the
spray $S$ is clear from context, we refer simply to geodesics.  Also, we
shall frequently consider noncompact manifolds because no spray can be
disprisoning on a compact manifold.  However, Corollary~\ref{qcov} may be
used to obtain results about compact manifolds for which the universal
covering is noncompact.

We refer to \cite{DRP1} for motivation, further general results, and
results specific to homogeneous sprays, and to \cite{DRP2} for more
examples.  Note that the sprays in \cite{DRP1} were positively 
homogeneous; the extension of those results to complete homogeneity is
straightforward, once the definition of homogeneous spray there is 
corrected to the one here.

\section{Exponential maps}\label{exp}
Let $S$ be a spray on $M$. We define the generalized exponential map{\em 
s\/} (plural!) $\exp^\ve$ of $S$ as follows.

First let $p\in M$, $v\in T_p M$, and $c$ be the 
unique $S$-geodesic such that
\begin{eqnarray*}
\ddc &=& S(\dc)\\
c(0) &=& p\\
\dc(0) &=& v
\end{eqnarray*}
Define
$$\exp^\ve_p(v) = c(\ve)$$
for all $v\in T_p M$ for which this makes sense.  From the existence of
flows ({\em e.\,g.,} \cite[p.\,175]{HS}), it follows that this is well
defined for all $\ve$ in some open interval $(-\ve_p,\ve_p)$, which in
general depends on $p$, and for all $v$ in some open neighborhood $U_p$ of
$0\in T_p M$, which in general depends on the choice of
$\ve\in(-\ve_p,\ve_p)$.  This defines $\exp^\ve_p$ at each $p\in M$.

Next, choose a smooth function $\ve :  M \to \R$ such that $\ve(p) \in
(-\ve_p,\ve_p)$ for every $p \in M$.  (The smoothness of $\ve$ is for our
later convenience: we want $\exp^\ve_p$ to be smooth in $\ve$ as well as 
in all other parameters.)  Then the global map $\exp^\ve$ is defined
pointwise by $(\exp^\ve)_p = \exp^{\ve(p)}_p$.  The domain of $\exp^\ve$
is a tubular neighborhood of the 0-section in $TM$ and the graph of $\ve$
lies in a tubular neighborhood of the 0-section in the trivial line bundle
$\R \times M$.

We have an example, given to us by J. Hebda, to show that it is possible 
that $\ve_p<1$ for every open neighborhood of $0\in T_p M$ if the spray is
inhomogeneous.
\begin{example}
Consider the spray on $\R$ given by
$$\ddx(t) = \pi\left(1 + \dx(t)^2 \right).$$
To integrate, we rewrite this as
$$\frac{d\dx}{1+\dx^2} = \pi\, dt$$
and obtain
$$\arctan\dx = \pi\, t+C_1\,.$$
Thus
$$\dx(t) = \tan \left(\pi\, t + C_1\right), \quad\dx(0) = \tan C_1$$
so
$$x(t) = \log\bigl|\sec\left(\pi\, t + C_1\right)\bigr| + C_2\,.$$
For $C_1>0$, $x$ cannot be continued beyond
\begin{eqnarray*}
\pi t + C_1 &=& \frac{\pi}{2}\,,\\[4pt]
t &=& \frac{1}{2} - \frac{C_1}{\pi} < 1\,.
\end{eqnarray*}
Therefore the usual exponential map of this spray is not defined ({\em 
i.\,e.,} at $t=\ve =1$) for any $C_1>0$. 
\end{example}

The closer the graph of $\ve$ gets to the 0-section of $\R \times M$, the 
larger the tubular neighborhood of the 0-section in $TM$ gets.
\begin{proposition}
For $\ve_1 \le \ve_2$, we have $\dom(\exp^{\ve_{\!\!\;2}}) 
\supseteq \dom(\exp^{\ve_{\!\!\;1}})$, attaining all of $TM$ for $\ve = 0$
when $\exp^0 = \pi$.\eop
\label{exp0}
\end{proposition}
This puts the bundle projection $TM\surj M$ in the interesting position of 
being a member of a one-parameter family of maps, all of whose other members
are local diffeomorphisms. (This is reminiscent of singular perturbations.)
\begin{theorem}
For every\/ $\ve$ such that\/ $0<|\ve|<\ve_p$, the generalized exponential 
map\/ $\exp^\ve_p$ is a diffeomorphism of an open neighborhood of\/ $0 \in
T_p M$ with an open neighborhood of\/ $p \in M$.
\end{theorem}
\begin{proof}
This follows from the flow theorems in ODE ({\em e.\,g.,} \cite[pp.\,175,
302]{HS}).
\end{proof}
Note that for $v \in T_p M$, $\exp^\ve_p v = \pi\P(\ve,v)$ where $\P$ is 
the local flow of $S$.

\begin{figure}
\begin{center}
\leavevmode
\def\epsfsize#1#2{.5#1}
\epsffile{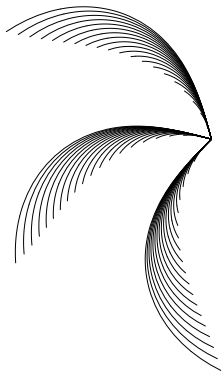}
\caption{\small curves $\exp^\ve_p(av)$ --- Each black curve is a geodesic 
with $0 < \ve < 3$ and $a$ and $v$ fixed.  From shortest to longest in each
plume, $a$ steps in increments of 0.05 from 0.05 to 1. In each plume, $v$ is
constant.  There are three implicit $a$-parameter curves readily located,
one along the endpoints of each of the three plumes.}\label{j1}
\end{center}
\end{figure}

\begin{figure}
\begin{center}
\leavevmode
\def\epsfsize#1#2{.6#1}
\epsffile{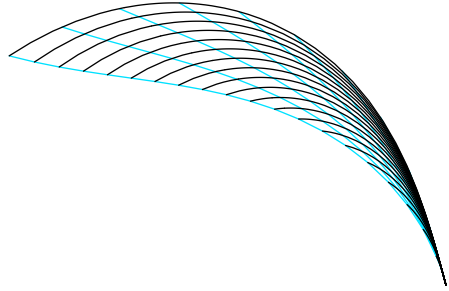}
\caption{\small curves $\exp^\ve_p(av)$ --- This is one plume from Figure 
\protect\ref{j1}. Each black curve is a geodesic and each gray curve is an
$a$-parameter curve. The new Jacobi fields are {\em along\/} the black 
curves but {\em tangent\/} to the gray curves.}\label{j2}
\end{center}
\end{figure}

For reference, we record the following obvious result.
\begin{lemma}
$\ve$ is a geodesic parameter; {\em i.\,e.,} the curve obtained by fixing 
$v$ and varying $\ve$ is a geodesic through $p$.\eop
\end{lemma}
Now consider another parameter $a$ as in
$$\exp^\ve_p(av)\,.$$
In general, $a$ will {\em not\/} be a geodesic parameter; {\em i.\,e.,} the
curve obtained by fixing $\ve$ and $v$ and varying $a$ is {\em not\/} a
geodesic through $p$.  See Figures \ref{j1} and \ref{j2} for a comparison.
\begin{proposition}
If $S$ is homogeneous, then $a$ as above is a geodesic parameter.
\end{proposition}
\begin{proof}
When $S$ is homogeneous, we can take $\ve=1$ and recover the usual 
exponential map, and then $a$ is the usual geodesic parameter.
\end{proof}
The $a$-parameter curves are interesting:  they are the integral curves for
our new Jacobi vector fields. These were mentioned in \cite{DRP2} and will
be studied in more detail elsewhere.  For now, we have the following
example.
\begin{example}
In $\R^2$, consider the spray given by $\S^i(x,y) = x^i$ for $i=1,2$. The 
geodesics are easily found to be $c(t) = v e^t + p$ where $v$ is the 
initial velocity and $p$ is the initial position.  We can use the usual 
exponential map since these curves are always defined for $t=1$. Thus we
obtain $\exp_{p}(v) = c(1) = v\,e + p$, regarding both $v$ and $p$ as
vectors in $\R^2$.

For the $a$-curves, we have $\exp_p(a\,v) = av\,e + p$, showing the 
difference between the two types quite clearly: the geodesics have 
exponential growth in velocity, while the $a$-curves have only linear
growth.
\end{example}
Finally, note that we could just as well define exponential-like maps
based on the $a$-curves and they would share most of the properties
of our new exponential maps.

\section{Connections and sprays}\label{cs}
In general, a {\em connection\/} on a manifold $M$ is a subbundle $\H$ of
the second tangent bundle $\pi_T:TTM\surj TM$ which is complementary to
the vertical bundle $\V = \ker(\pi_* :TTM\surj TM)$, so
\begin{equation}
TTM = \H\ds\V\, .
\label{ceq}
\end{equation}
We note there are two vector bundle structures on $TTM$ over $TM$, denoted
here by $\pi_T$ and $\pi_*$.  While $\V$ is always a subbundle with
respect to both \cite[pp.\,18,20]{P}, $\H$ is a subbundle with respect to
$\pi_*$ if and only if the connection is linear \cite[p.\,32]{B}.

Recall that quadratic sprays correspond to linear connections.  In terms of
the horizontal bundle $\H\! ,$ linearity is expressed as
$$\H_{av} = a_* \H_v $$
for $a\in\R$ and $v\in TM$. Thus one has 
\begin{equation}
\H_{av} = a_*a^{m-1}\H_v
\label{hceq}
\end{equation}
as the second defining equation, together with (\ref{ceq}), of a
connection that is $h(m)$. 

Note that for an $h(m)$ semispray $S$ with integral $m$, Grifone's
\cite{G1} associated (generalized) Chris\-toffel symbols $\Gamma$ are
$h(m-1)$, appropriately.  See (\ref{ours}) below for our version, which
allows for nonlinear, including inhomogeneous, connections.

Here is the spray induced by a connection.  We shall call it the {\em
geodesic spray\/} associated to the connection and its geodesics the {\em
geodesics\/} of the connection.
\begin{theorem}\label{cis}
For each connection $\H\! ,$ there is an induced spray $S$ given by
$$S(v) = \pi_*\big|^{-1}_{\H_v} (v)\, ,$$
where $\pi:TM\surj M$ is the natural projection and $v\in TM$. We write 
$\H\vdash S$ to denote this relationship.
\end{theorem}
\begin{proof}
As in the first paragraph of Poor's proof of 2.93 \cite[p.\,95]{P}, it is
easily verified that $S$ so defined is a spray.  Indeed, $S$ is a section of
$\pi_*$ by construction, and $S$ is a section of $\pi_T$ because $\H$ is a
subbundle with respect to $\pi_T$.
\end{proof}
It is clear that this spray is horizontal, so compatible with the given 
connection.

Unfortunately, when the connection is $h(m-1)$ this spray is not
homogeneous {\em as a spray;} it is only an $h(m)$ {\em vector field\/} on
$TM$.  In order to avoid this problem, we must consider a new type of
(partial) homogeneity for connections.
\begin{definition}\label{dvhc}
A connection $\H$ on $TM$ is\/ {\em vertically homogeneous} of degree $m$, 
denoted by $vh(m)$, if and only if
\begin{equation}\label{vhceq}
\H_{av} = a_* \av^{m-1} \H_v
\end{equation}
where $\av^n$ denotes scalar multiplication by $a^n$ in the vertical 
component and $a$ in the horizontal component in some (hence any)
admissible local coordinates.
\end{definition}
More explicitly, $a_* \av^{m-1} (x,y,X,Y) = (x,ay,aX,a^m Y)$ in admissible 
coordinates. Note that $h(m)$ and $vh(m)$ coincide only for $m=1$, the 
linear connections.
\begin{proposition}\label{hgs}
If $\H$ is a connection with geodesic spray $S$, then $S$ is $h(m)$ if 
and only if $\H$ is $vh(m-1)$.
\end{proposition}
\begin{proof}
That $S$ is $h(m)$ if $\H$ is $vh(m-1)$ follows as in the second
paragraph of Poor's proof of 2.93 \cite[p.\,95]{P}, {\em mutatis 
mutandis;} the converse results from a similar calculation.
\end{proof}

Connections may also be seen as sections of the bundle $G_H(TTM)$ of all
possible horizontal spaces, a subbundle of the Grassmannian bundle
$G_n(TTM)$.  To see what structure $G_H(TTM)$ has, consider $\R^{2n} =
\R^n\ds\R^n$ as the model fiber of $TTM$ and regard the first summand as
horizontal, the second as vertical.  With $GL_{2n}$ as the structure group
of $TTM$, we want the subgroup $A_H$ that preserves the vertical space and
maps any one horizontal space into another.  This can be conceived as
occurring in two steps.  First, we may apply any automorphisms of the
vertical and horizontal spaces separately.  Second, we may add vertical
components to horizontal vectors to obtain the new horizontal space.
$$ \arraycolsep=.5em \left[ 
\begin{array}{cc} I & 0 \\ \gl_n & I \end{array} \right] \cdot \left[
\begin{array}{cc} GL_n & 0 \\ 0 & GL_n \end{array} \right] $$
Our group $A_H$ is thus found to be a semidirect product entirely
analogous to an affine group.  The action is transitive and the right-hand
factor is the isotropy group of any fixed horizontal space, so the model
fiber for $G_H(TTM)$ is the resulting homogeneous space.  The induced
operation on representatives being given by
$$ \arraycolsep=.5em \left[ 
\begin{array}{cc} I & 0 \\ A & I \end{array} \right] \cdot \left[
\begin{array}{cc} I & 0 \\ B & I \end{array} \right]  = \left[
\begin{array}{cc} I & 0 \\ A+B & I \end{array} \right], $$
it follows that $G_H(TTM)$ is an affine bundle (bundle of affine spaces, 
{\em vs.} vector spaces).  Thus a connection, being a section of this 
bundle, provides a choice of distinguished point in each fiber, hence a 
vector bundle structure on this affine bundle.

If we wish to consider only those connections compatible with a given
spray, we just replace arbitrary elements of $\gl_n$ with those having a
first column comprised entirely of zeros.  Note that this yields an affine
subbundle $G_H^S(TTM)$ of $G_H(TTM)$, with fibers being pencils of 
possible horizontal spaces.
\begin{theorem}\label{cce}
Given a spray $S$ on $M$, there exists a compatible connection $\H$ in 
$TTM$.
\end{theorem}
Since the fibers of $G_H^S(TTM)$ are contractible, this is an easy 
exercise in obstruction theory \cite[Ch.\,8]{DP}; however, an explicit 
construction is desirable.
\begin{proof}
Let $\P$ denote the local flow of $S$ and $\g$ an integral curve of $S$
with $\g(0) = v \in T_p M$.  The basic idea is to use $S$ and $\P$ to {\em
define\/} notions of {\em horizontal\/} and {\em parallel\/} which will
coincide with the usual ones along $\g$ for any $\H \nd S$.  This is
essentially the same as the usual construction \cite{P}.  The problem is
that for inhomogeneous $S$, the ray $\{tv\}$ in $T_p M$ does {\em not\/}
exponentiate to a geodesic in $M$.

To remedy this, we proceed as follows. For each $v\in T_p M$, choose 
$\ve_v$ so that $\exp_p^{\ve_{\!\!\:v}}v$ is defined. Such $\ve_v$ exist 
by Proposition \ref{exp0}. For $0 \le t \le \ve_v$, define 
\begin{equation}
\alpha_v(t) = \left(\exp_p^{\ve_{\!\!\:v}}\right)^{-1}\exp_p^t v\,.
\end{equation}
Then $\alpha_v(0) = 0$, $\alpha_v(\ve_v) = v \in T_p M$, and $\alpha_v$ 
exponentiates to the geodesic with initial condition $v$ at $p$. Note that
if $S$ is homogeneous, then $\alpha_v(t) = tv$.

We have a vector bundle map $\cJ:\pi^*TM\to \V$ which is an isomorphism on
fibers.  It is one version of canonical parallel translation on a vector
space, identifying the tangent space at each point with the vector space
itself.  Now, for each $w \in T_p M$ define
\begin{equation}
\H_w = \left\{ \left.\frac{d}{dt}\right|_{t=0} \pi_* 
\P_{t*}\,\cJ_{\alpha_{\!\!\;v}(t)}w \Bigm| v \in T_p M \right\} .
\end{equation}
Clearly, this does not depend on the choices of $\ve_v$ made earlier.
(Note we are evaluating at 0.)  If $S$ is quadratic, it is easy to check
that this coincides with the usual construction as found in
\cite[pp.\,96--97]{P}, since in that case $\exp_p tv = \pi\P(t,v)$ for $v
\in T_p M$.  The proof that $\H$ so defined is a connection and that $\H
\nd S$ follows Poor's proof of 2.98 \cite[pp.\,97--99]{P} {\em mutatis
mutandis.}
\end{proof}
These connections will be our ``standard"---our generalization of 
torsion-free linear connections; {\em cf.} equation (\ref{new}), 
Definition \ref{torsfree} and after.

We further note that admissible atlases correspond to certain reductions 
of the structure group of $TTM$ from $GL_{2n}$ to $GL_n \ds GL_n$, those 
which in turn correspond to direct-sum decompositions of $TTM$ in which 
one of the summands is the vertical bundle $\V$ (and the other is perforce
a horizontal bundle), hence to connections in $TTM$. Thus any homogeneous 
spray $S$ comes with a particular associated compatible connection $\H$,
the one corresponding to the associated admissible atlas; {\em cf.} after
Definition \ref{hms}.  Note, however, that it may not be the one naturally
associated by the preceding construction.

Here is an alternative, axiomatic characterization of a connection in terms 
of the horizontal projection $H$.
\begin{itemize}
\item[\bf C1 ] $H$ is a smooth section of $\End(TTM)$ over $TM$.

\item[\bf C2 ] $H^2 = H$.

\item[\bf C3 ] $\ker H = \V$.
\end{itemize}
Then $\H = \im H$ is the horizontal bundle. Vertical homogeneity is 
expressed with an optional axiom.
\begin{itemize}
\item[\bf Ch ] $H$ is $vh(m)$ if and only if $H_{av}a_* = a_*\av^{m-1}
H_v$ for all $v\in TM$ and $a\in\R$ ($v\in TM-0$ and $a\ne 0$ for $m<1$).
\end{itemize}
Homogeneous connections may be similarly axiomatized.

There is another natural vector bundle map $K:\V\to TM$ respecting $\pi_T$
which is an isomorphism on fibers, another version of canonical
parallel translation of a vector space.  Using this, we define a
connection map or connector for an arbitrary connection and thence a
covariant derivative.
\begin{definition}\label{kap}
For a connection $\H\! ,$ define the associated\/ {\em connector} $\kappa 
:TTM\to TM : z\mapsto K(z - H_v z)$ for $z\in T_v TM$.
\end{definition}
\begin{proposition}\label{kvbm}
The connector $\kappa$ is a vector bundle map respecting $\pi_T$ but\/ {\em
not} $\pi_*$ in general.  It respects $\pi_*$ if and only if the connection
is linear.
\end{proposition}
\begin{proof}
As in Poor \cite[p.\,72f\,]{P}, {\em mutatis mutandis.}
\end{proof}
According to Besse \cite[p.\,32f\,]{B}, a {\em symmetric\/} connector
(connection) is invariant under the natural involution $J$ of $TTM$.
Clearly this is possible only for linear connections.

Now we are ready for the main event.
\begin{definition}\label{cd}
The\/ {\em covariant derivative} associated to the connection\/ $\H$ is 
the operator defined by 
$$\del{u}v = \kappa(v_* u) = K(v_* u - H_v v_* u)$$
and is tensorial in $u$ but\/ {\em nonlinear} (in general) in $v$.
\end{definition}
This last comes from the general lack of respect for the $\pi_*$ structure
by $\H\! ,$ $H\!$, and $\kappa$.
\begin{example}\label{hex}
We always have $\del{0}v = 0$.  For any $vh(m)$ connection, $\del{u}av =
K(a_*v_*u - H_{av}a_*v_*u) = aK(v_*u - \av^{m-1}H_v v_*u)$, and similarly
for homogeneous ones. So (vertically) homogeneous connections do not
differ significantly from linear ones.  In particular, $\del{u}0 = 0$ for
all $u$ for all (vertically) homogeneous connections; in fact, they all
have the same horizontal spaces along the 0-section of $TM$, namely the
subspaces tangent to it ({\em i.\,e.,} those in the image of $0_* : TM\to
TTM$).  We call all such connections sharing this property {\em
0-preserving;} they differ minimally from (vertically) homogeneous
(including linear) connections.  In contrast, connections with $\del{u}0
\ne 0$ for even some $u$ are much farther from linear; we call them {\em
strongly nonlinear.} See Figure \ref{sch} for a schematic view.
\end{example}
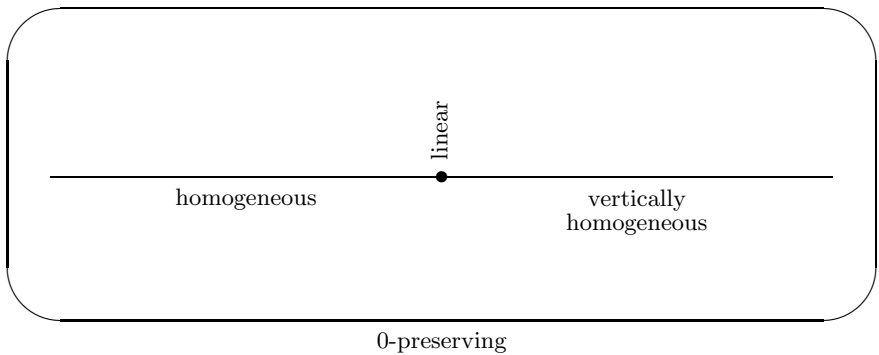
\begin{figure}
\begin{center}
\setlength{\unitlength}{.15em}
\begin{picture}(200,75)
\put(10,50){\line(1,0){180}}
\put(100,50){\makebox(0,0){$\bullet$}}
\put(55,45){\makebox(0,0){\footnotesize homogeneous}}
\put(145,42){\makebox(0,0){$\mbox{\footnotesize vertically}\atop
   \mbox{\footnotesize homogeneous}$}}
\put(97.5,54){\psrotate{\mbox{\footnotesize linear}}}
\put(100,53){\oval(200,72)}
\put(100,12){\makebox(0,0){\footnotesize 0-preserving}}
\end{picture}\end{center}
\caption{\small Each set of connections is closed with empty interior in 
the next:  linear in homogeneous, linear in vertically homogeneous, linear
and homogeneous in 0-preserving, linear and vertically homogeneous in
0-preserving, linear and homogeneous and vertically homogeneous in
0-preserving, 0-preserving in the whole.  The strongly nonlinear
connections may be visualized as a 3-d cloud containing the 0-preserving
ones.}\label{sch}
\end{figure}

As usual, $\frak{X}$ denotes the vector fields on $M$.  We recall the
vector bundle map $\cJ:\pi^*TM\to \V$ which is an isomorphism on fibers
and a version of canonical parallel translation on a vector space.
\begin{theorem}\label{corresp}
There is a bijective correspondence between (possibly nonlinear) 
connections $\H$ and our (possibly nonlinear) covariant derivatives $\D$ 
on $TM$.
\end{theorem}
\begin{proof}
It suffices to show that we can reconstruct $\H$ from its associated 
covariant derivative $\D$.  For each $u\in T_p M$, define
$$ \bar{\H}_u = \{U_*v - \cJ_u\del{v}U\mid U\in\frak{X}, U_p = u, v\in 
T_pM\} $$
and form the subbundle $\bar\H$ in $TTM$ in the obvious way.  It is easy
to see that $\bar\H$ is complementary to $\V$ as required, hence a
connection.  That $\bar\H$ is smooth is straightforward.  Finally, $\bar\H
= \H$ from this construction and the construction of $\D$ from $\H\!$.
\end{proof}
Compare \cite[p.\,77, proof of 2.58]{P}.  Thus as usual, we may refer
indifferently to $\H$ or its associated $\D$ as the connection.

Generalized Christoffel symbols may be introduced through
\begin{equation}
\left(KH_v v_*u\right)^k = \Gamma^k_i(v)u^i\,,\label{ours}
\end{equation}
making manifest the tensoriality in $u$. Here are some examples of their 
use. Observe that $(Kv_* u)^k = u^i\partial_i v^k$ so that
\begin{equation}
\left(\del{u}v \right )^k =  u^i\partial_i v^k - \Gamma^k_i(v)u^i
\label{lccd}
\end{equation}
is the covariant derivative. The geodesic equation is
\begin{equation}
\ddot{c}^{\,k} = \Gamma^k_i(\dot{c})\dot{c}^i ,
\label{lcgeq}
\end{equation}
which means that
\begin{equation}
\S^k(\dot{c}) = \Gamma^k_i(\dot{c})\dot{c}^i
\label{lcsg}
\end{equation}
for $S$ the spray induced by the connection $\D$. Note that we can write 
{\em any\/} $\S^k(v)$ in the form 
\begin{equation}
\label{new}
\S^k(v) = \Gamma^k_i(v)v^i\, ,
\end{equation}
although $\Gamma$ may be less well-behaved than $\S$.  In this way we can
obtain the standard (``torsion-free'') connection associated to $S$ by our
generalized APS construction (proof of Theorem \ref{cce}).

We obtain the usual relation between two notions of geodesic.
\begin{theorem}\label{cg=sg}
A curve $c$ is a geodesic of\/ $\H$ if and only if\/ $\del{\dc}\,\dc = 0$.
\end{theorem}
\begin{proof}
$\del{\dc}\,\dc = \kappa(\dc_*\dc) = K(\dc_*\dc - H_{\dc}\,\dc_*\dc) =
K(\dc_*\dc - S(\dc))$ by the construction of $S$ in Proposition \ref{cis}.
Now all we have to do is identify $\dc_*\dc$ as $\ddc$ and recall that $K$
is an isomorphism on fibers.
\end{proof}

Curvature is readily handled.  Let $\H$ be a connection on $M$.  The {\em
horizontal lift\/} of a vector field $U$ on $M$ is defined as usual and
denoted by $\bU$.
\begin{definition}\label{curv}
Given vector fields\/ $U$ and\/ $V$ on $M$, the\/ {\em curvature operator} 
$R(U,V):TM\to TM$ is defined by
$$R(U,V)w = \kappa\left([\bV,\bU]_w\right)$$
for all $w\in TM$.  It is tensorial in the first two arguments, but\/
{\em nonlinear} (in general) in the third.
\end{definition} 
The arguments are reversed on the right in order to obtain the usual formula
in terms of the associated covariant derivative,
$$R(U,V)W = \del{U}\del{V}W - \del{V}\del{U}W - \del{[U,V]}W\, ,$$
as one may verify readily.  It is also easy to check that this curvature 
vanishes if and only if $\H$ is integrable, thus justifying our definition.
\label{more?}

Torsion is considerably more obscure.  Consider two (possibly nonlinear)
connections $\bar\H$ and $\H$ on $TM$ with corresponding (possibly
nonlinear) covariant derivatives $\bar\D$ and $\D$.
\begin{definition}\label{dif}
Given two covariant derivatives $\bar\D$ and $\D$, define the\/ {\em 
difference operator} $\cD = \bar\D - \D$.
\end{definition}
We think of $\cD$ as having two arguments, $\cD(u,v) = \bdel{u}v - 
\del{u}v$.  It is always tensorial in $u$, but is {\em nonlinear\/} (in
general) in $v$.  Alternative notations include $\cD_u v$ and $\cD(u)v$
when one wishes to emphasize certain aspects.

We define the {\em covariant differential\/} as usual {\em via\/} $(\D v)u
= \del{u}v$. As an operator, $\D v$ is still linear in its argument $u$.
\begin{proposition}\label{bar}
For all $v\in TM$, $\bar{\H}_v = \{z - \cJ_v\cD(\pi_* z,v)\mid 
z\in\H_v\}$.
\end{proposition}
\begin{proof}
Let $v\in T_p M$, $z\in\H_v$, $V\in\frak{X}$ such that $(\D V)_p = 0$ and 
$V_p = v$. Thus if $u = \pi_*z\in T_p M$, then $z=V_*u\in\H_v$. Now
$$ \bar\kappa V_*u = \bdel{u}V = \del{u}V + \cD(u,v) = \cD(u,v) = 
\bar\kappa\cJ_v\cD(u,v)\,,$$
so $\bar\kappa\left( z - \cJ_v\cD(u,v)\right) = 0$ and $z - \cJ_v\cD(u,v) 
\in \bar{\H}_v$.

Since $\pi_*$ is an isomorphism of the horizontal spaces $\bar\H_v$ and
$\H_v$ with $T_p M$ and $\pi_*z = \pi_*\left( z - \cJ_v\cD(u,v)\right)$,
this yields all of $\bar\H_v$.
\end{proof}
\begin{proposition}\label{=geos}
Two connections on $TM$ have the same geodesic spray if and only if their 
associated difference operator is alternating (vanishes on the diagonal 
of\/ $TM\ds TM$).
\end{proposition}
\begin{proof}
For each $v\in TM$, $S_v = \pi_*\big|^{-1}_{\H_v} (v)$ while $\bar{S}_v
= \pi_*\big|^{-1}_{\bar{\H}_v} (v) = \pi_*\big|^{-1}_{\H_v} (v) - 
\cJ_v\cD(v,v)$. Therefore $\bar{S} = S$ if and only if $\cD(v,v) = 0$ for 
all $v\in TM$.
\end{proof}
For {\em linear\/} connections, $\cD$ is {\em bilinear\/} and alternating
is equivalent to antisymmetric (or, skewsymmetric).  In general, of
course, this does not hold.

The familiar formula for torsion $T(u,v) = \del{u}v - \del{v}u - [u,v]$ is
not linear (let alone tensorial) in either argument.  Thus the usual trick
to get a torsion-free linear connection, replacing $\D$ by $\bar\D = \D -
\half T$, will not work for our nonlinear connections.  Indeed, $\bar\D$
and $\D$ seem to have the same geodesics and $\bar\D$ is formally
torsion-free, {\em but\/} the new $\bar\D$ is {\em not\/} one of our
nonlinear covariant derivatives:  $\bdel{u}v$ is {\em not\/} tensorial in
$u$.

\begin{definition}\label{op}
Denote by\/ $\Op TM$ the smooth maps $TM\to TM$ which preserve fibers
(\,{\em i.\,e.,} commute with the projection onto $M$).  We shall write\/ 
$\Op_m TM$ when the maps are $h(m)$ on each fiber.
\end{definition}
Note they are smooth on fibers, but not necessarily linear.  When they are
fiberwise linear, we have $\End TM$ as in the usual formulation.  As is
$\End TM$, $\Op TM$ is a vector bundle over $M$ and a Lie algebra.
(Addition, scalar multiplication, and composition---thus commutators---are
well defined and preserve fibers.) By analogy to the linear theory, we 
usually think of a covariant derivative $\D$ as a section (over $M$) of 
$\Hom\left(TM,\Op TM\right)$ and $\cD$ as an $\Op TM$-valued 1-form on 
$M$.

A replacement $\cT$ for torsion must also be an element of $A^1(M,\Op
TM)$ in order for it to play the same role in general that torsion does
for linear connections.  For then, given any $\cT\in A^1(M,\Op TM)$,
$\bar\D = \D + \cT$ is another nonlinear covariant derivative of our type.
We want to choose $\cT$ so that $\bar\D$ has the same geodesics as $\D$ 
but is as analogous to a torsion-free linear connection as possible.

As we noted immediately after the proof of Theorem \ref{cce}, what we 
shall do is one of the classic mathematical gambits: turn a theorem into
a definition.
\begin{definition}
We define the connections constructed in the proof of Theorem \ref{cce}
to be the (generalized)\/ {\em torsion-free} connections.
\label{torsfree}
\end{definition}
(We refer to Poor \cite[pp.\,101--102]{P} for the relation to the classic 
Ambrose-Palais-Singer correspondence.)  Equivalently, we are regarding the
usual torsion formula as derived from the difference operator (difference
tensor in the linear case) construction; {\em cf.} \cite[pp.\,99--100]{P}.

Now we may construct the (generalized) torsion of a (possibly nonlinear) 
connection $\H$ with corresponding (possibly nonlinear) covariant 
derivative $\D$. By Theorem \ref{cis}, $\H$ induces a (unique, 
generalized) spray $S$. Use the proof of Theorem \ref{cce} to construct 
the connection $\hat\H$ from $S$. By Theorem \ref{corresp} there is a 
unique covariant derivative $\hat\D$ corresponding to $\hat\H$. Let $\cD =
\D - \hat\D$ be the difference operator, so $\hat\D = \D - \cD$ is 
(generalized) torsion-free.
\begin{definition}
Using the preceding notations, the (generalized)\/ {\em torsion} of\/ $\D$
is defined by $\cT = 2\cD = 2\left(\D - \hat\D\right)$.
\label{tors}
\end{definition}
The factor of two here and the subtraction order make verification that 
this reduces to classical torsion in the linear case immediate, and 
preserves the traditional formula $\hat\D = \D - \half\cT$ for the 
associated torsion-free connection.

\section{Geodesic connectivity and stability}\label{gcs}
In \cite{DRP1}, we defined a spray to be {\sc LD} if and only if its usual
exponential map is a local diffeomorphism.  For some results there, we used
the fact that the geodesics of such sprays give normal starlike
neighborhoods of each point in $M$.  (In fact, the $a$-curves also give 
such neighborhoods, as is easily seen.) These results now immediately extend 
to all sprays.  For convenience, we state them here.
\begin{proposition}
Let $M$ be a manifold with a pseudoconvex and disprisoning spray $S$. If\/ 
$S$ has no conjugate points, then $M$ is geodesically connected.
\label{p5}
\end{proposition}

Let $M$ be a manifold with a spray $S$ and let $\widetilde {M}$ be a
covering manifold.  If $\phi:\widetilde{M}\rightarrow M$ is the covering
map, then it is a local diffeomorphism.  Thus $\tilde{S} =
(\phi_{\ast})^{\ast}S$ is the unique spray on $\widetilde{M}$ which covers
$S$, geodesics of $\tilde{S}$ project to geodesics of $S$, and geodesics of
$S$ lift to geodesics of $\tilde{S}$.  Also, $S$ has no conjugate points if
and only if $\tilde{S}$ has none.  The fundamental group is simpler, and
$\tilde{S}$ may be both pseudoconvex and disprisoning even if $S$ is
neither.
\begin{corollary}
Let $M$ be a manifold with a pseudoconvex and disprisoning spray $S$ and let
$\widetilde{M}$ be a covering manifold with covering spray $\tilde{S}$.  If
$\tilde S$ has no conjugate points, then both $\widetilde M$ and $M$ are
geodesically connected.  \label{qcov}
\end{corollary}
\begin{theorem}
Let $S$ be a pseudoconvex and disprisoning spray on $M$.  If $S$ has no
conjugate points, then for each $p \in M$ the exponential maps of $S$ at $p$
are diffeomorphisms.
\end{theorem}
We remark that none of these results require (geodesic) completeness of
the spray $S$.

We now consider the joint stability of pseudoconvexity and disprisonment
for (general) sprays in the fine topology.  Because each linear connection
determines a (quadratic) spray, Examples 2.1 and 2.2 of \cite{BP4} show
that neither condition is separately stable.  (Although \cite{BP4} is
written in terms of principal symbols of pseudodifferential operators, the
cited examples are actually metric tensors).  We shall obtain $C^{0}$-fine
stability, rather than $C^{1}$-fine stability as in \cite {BP4}, due to
our effective shift from potentials to fields as the basic objects.  The
proof requires some modifications of that in \cite {BP4}; we shall
concentrate on the changes here and refer to \cite {BP4} for an outline
and additional details.

Rather than considering $r$-jets of functions, we now take $r$-jets of
sections in defining the Whitney or $C^r$-fine topology as in Section 2 of
\cite{BP4}.  Let $h$ be an auxiliary complete Riemannian metric on $M$.
Thus we look at the $C^r$-fine topology on the sections of $TTM$ over $TM$.

If $\g_1$ and $\g_2$ are two integral curves of a spray $S$ with $\g_1(0) =
(x,v)$ and $\g_2(0) = (x,\l v)$ for some positive constant $\l$, then the
inextendible geodesics $\p\circ\g_1$ and $\p\circ\g_2$ no longer differ only
by a reparametrization.  Thus, in contrast to \cite{BP4}, we must now
consider an integral curve for each non-zero tangent vector at each point of
$M$.  Note this also means that we can no longer use the $h$-unit sphere 
bundle to obtain compact sets in $TM$ covering compact sets in $M$.

Observe that the equations of geodesics involve no derivatives of $S$.
Thus if $\g :[0,a]\rightarrow TM$ is a fixed integral curve of $S$ in $TM$
with $\g(0) = v_0 \in TM$ and if $\g' :[0,a]\rightarrow TM$ is an integral
curve of $S'$ in $TM$ with $\g'(0) = v$, then $d_h \left( \p\circ\g(t),
\p\circ\g'(t)\right) < 1$ for $0\le t\le a$ provided that $v$ is
sufficiently close to $v_0$ and $S'$ is sufficiently close to $S$ in the
$C^0$-fine topology.  This and the $\sigma$-compactness of $T K_1$ when 
$K_1$ is compact yield the following result.
\begin{lemma}
Assume $K_{1}$ is a compact set contained in the interior of the compact
set $K_{2}$, $V$ is an open neighborhood of $K_2$, $S$ is a disprisoning
spray, and let $\epsilon > 0$.  There exist countable sets $\{v_i\}
\subseteq TK_{1}$ of tangent vectors and $\{\delta_{i}\}$ and $\{a_{i}\}$
of positive constants such that if $S'$ is in a $C^{0}$-fine
$\epsilon$-neighborhood of $S$ over $V$, then the following hold:
\begin{enumerate}
\item if $c$ is an inextendible $S$-geodesic with $c(0)$ in a
$\delta_{i}$-neighborhood of $v_{i}$, then $c[0,a_{i}] \subset V$ and $c
(a_{i}) \in V-K_{2}$;

\item If $c^{\prime}$ is an inextendible $S^{\prime}$-geodesic with
$\dot{c}^{\prime}(0)$ in a $\delta _{i}$-neighborhood if $v_{i}$, then
$c^{\prime}[0,a_{i}] \subset V $ and $c^{\prime}(a_{i}) \in V-K_{2}$;

\item Two inextendible geodesics, $c$ of $S$ and $c^{\prime}$ of
$S^{\prime}$ with $\dot{c}(0)$ and $\dot{c}'(0)$ in a $\delta
_{i}$-neighborhood of $v_{i}$, remain uniformly close together for $0\leq t
\leq a_{i}$;

\item The union of all the $\delta_{i}$-neighborhoods of the $v_{i}$ 
covers $TK_{1}$.\eop
\end{enumerate} 
\end{lemma}

Continuing to follow \cite{BP4}, we construct the increasing sequence of
compact sets $\{ A_n\}$ which exhausts $M$ and the monotonically
nonincreasing sequence of positive constants $\{ \epsilon_n\}$.  The only
additional changes from \cite[p.\,17f\,]{BP4} are to use integral curves of
$S$ in $TM$ instead of bicharacteristic strips in $T^*M$.  No other
additional changes are required for the proof of the next result either.
\begin{lemma}
Let $S$ be a pseudoconvex and disprisoning spray and let $S^{\prime}$ be
$\delta$-near to $S$ on $M$.  If $c^{\prime}:(a,b)\rightarrow M$ is an
inextendible $S'$-geodesic, then there do not exist values $a < t_{1} <
t_{2} < t_{3} < b$ with $c^{\prime} (t_{1}) \in A_{n}$, $c ^{\prime} (t_{3})
\in A_{n}$, and $c^{\prime} (t_{2}) \in A_{n+4} - A_{n+3}$.\eop \label{3.2}
\end{lemma}

Now we establish the stability of pseudoconvex and disprisoning sprays by
showing that the set of all sprays in $\spr(M)$ which is pseudoconvex and
disprisoning is an open set in the $C^0$-fine topology.  The only changes
needed from the proof of Theorem 3.3 in \cite[p.\,19]{BP4} are replacing
principal symbols by sprays, bicharacteristic strips by integral curves,
$S^*A_n$ by $TA_n$, and references to Lemma 3.2 there by references to
Lemma~\ref{3.2} here.
\begin{theorem}
If $S\in\spr(M)$ is a pseudoconvex and disprisoning spray, then there is
some $C^{0}$-fine neighborhood $W(S)$ in $\spr(M)$ such that each $S' \in
W(S)$ is both pseudoconvex and disprisoning.\eop \label{spd}
\end{theorem}
\begin{corollary}
If $M$ is a pseudoconvex and disprisoning pseudoriemannian manifold, then
any (possibly nonlinear) connection on $M$ which is sufficiently close to
the Levi-Civita connection is also pseudoconvex and disprisoning.\eop
\end{corollary}


\begin{thebibliography}{99}\frenchspacing

\bibitem{APS}
W. Ambrose, R.\,S. Palais and I.\,M. Singer, Sprays, {\it Anais Acad.
Brasil Ci\^ enc.} {\bf 32} (1960) 163--178.

\bibitem{AA}
P.L. Antonelli and M. Anastaseie, {\it The Differential Geometry of 
Lagrangians which Generate Sprays.} 
Dordrecht:  Kluwer, 1996.

\bibitem{AIM}
P.\,L. Antonelli, R.\,S. Ingarden and M.\,S. Matsumoto, {\it The Theory
of Sprays and Finsler Spaces with Applications in Physics and Biology.}
Dordrecht:  Kluwer, 1993.

\bibitem{A}
G.\,S. Asanov, {\it Finsler Geometry, Relativity and Gauge Theories.}
Dordrecht: Kluwer, 1985.

\bibitem{BCS}
D. Bao, S.-S. Chern, and Z. Shen, {\it An Introduction to Riemann-Finsler 
Geometry.} New York: Springer, 2000.

\bibitem{BP4}
J.\,K. Beem and P.\,E. Parker, Whitney Stability of Solvability, {\it Pac.
J. Math.} {\bf 116} (1985) 11--23.

\bibitem{BP6}
J.\,K. Beem and P.\,E. Parker, Pseudoconvexity and geodesic connectedness,
{\it Ann. Mat. Pura Appl.} {\bf 155} (1989) 137--142.

\bibitem{B}
A.\,L. Besse, {\it Manifolds all of Whose Geodesics are Closed}. New York:
Springer-Verlag, 1978.

\bibitem{BPo}
Misha Bialy and Leonid Polterovich, Geodesics of Hofer's metric
on the group of Hamiltonian diffeomorphims, {\it Duke Math J.} {\bf  76} 
(1994) 273--292.

\bibitem{BC}
F. Brickell and R.\,S. Clark, {\it Differentiable Manifolds}. New York: Van
Nostrand, 1970.

\bibitem{DR1}
L. Del Riego, 1-homogeneous sprays in Finsler manifolds, in {\it Global 
differential geometry: the mathematical legacy of Alfred Gray,} eds. 
Marisa Fern\'andez and Joseph A. Wolf. Contemp. Math. 288. Providence:
AMS, 2001. pp.\,411--414.

\bibitem{DRD}
L. Del Riego and C.\,T.\,J. Dodson, Sprays, universality and stability,
{\it Math. Proc. Camb. Phil. Soc.} {\bf 103} (1988) 515--534.

\bibitem{DRP1}
L. Del Riego and P.\,E. Parker, Pseudoconvex and disprisoning homogeneous
sprays, {\it Geom. Dedicata} {\bf 55} (1995) 211--220.

\bibitem{DRP2}
L. Del Riego and P.\,E.  Parker, Some nonlinear planar sprays, in {\it
Nonlinear Analysis in Geometry and Topology,} ed.  T.\,M.  Rassias.  Palm
Harbor: Hadronic Press, 2000.  pp.\,21--52.

\bibitem{DP}
C.\,T.\,J.  Dodson and P.\,E.  Parker, {\it A User's Guide to Algebraic
Topology.} Boston: Kluwer Academic Publishers, 1997.

\bibitem{D}
P. Dombrowski, On the geometry of the tangent bundle, {\it J. reine angew.
Math.} {\bf 210} (1962) 73--88.

\bibitem{G}
J. Grifone, Connexions non lin\'eaires conservatives, {\it C.\,R. Acad. Sci.
Paris S\'er. A Math.} {\bf 268} (1969) 43--45.

\bibitem{G1}
J. Grifone, {\it Structure Presque Tangent et Connexions non Homog\`enes.}
Th\`ese $3^{\mbox{\scriptsize\`eme}}$ cycle, Universit\'e de Grenoble, 
1971.

\bibitem{HS}
M.\,W. Hirsch and S. Smale, {\it Differential Equations, Dynamical Systems, 
and Linear Algebra.} New York: Academic Press, 1974.

\bibitem{KV}
J. Klein and A. Voutier, Formes ext\'erieures g\'eneratrices de sprays, {\it
Ann. Inst. Fourier} {\bf 18} (1968) 241--260.

\bibitem{LR}
M. de Le\'on and P. Rodr\'\i guez, {\it Methods of Differential Geometry in
Analytical Mechanics.} Amsterdam: North-Holland, 1989.

\bibitem{Po}
Leonid Polterovich, Geometry on the group of Hamiltonian diffeomorphisms,
{\it Doc. Math.} J. Extra Volume ICM II (1998) 401-410.

\bibitem{P}
W.\,A. Poor, {\it Differential Geometric Structures.} New York: McGraw-Hill,
1981.

\bibitem{R}
H. Reckziegel, Generalized sprays and the theorem of Ambrose-Palais-Singer,
in {\it Geometry and Topology of Submanifolds V,} ed.  F. Dillen, L.
Vrancken, L. Verstraelen, and I. Van de Woestijne.  River Edge:  World
Scientific, 1993, pp.\,242--248.


\bibitem{V}
A. Vondra, Sprays and Homogeneous Connections on $\R \times TM$, {\it Arch.
Math. (Brno)} {\bf 28} (1992) 163--173.

\end{thebibliography}
\end{document}